\documentclass[letterpaper, 10 pt, journal, twoside]{IEEEtran}
\usepackage{graphics} 
\usepackage{graphicx}
\usepackage{amsmath}
\usepackage{amssymb}
\usepackage{xcolor}
\usepackage{cite}
\usepackage{epstopdf}
\usepackage{tikz}
\usepackage[justification=centering]{caption}
\usetikzlibrary{shapes,arrows}

\title{A Generalization of the Classical Kelly Betting\\ Formula to the Case of Temporal Correlation}

\author{Joseph~D.~O'Brien, Kevin~Burke, Mark~E.~Burke, and B.~Ross~Barmish %
\thanks{Joseph D. O'Brien is a doctoral candidate
		in MACSI, Department of Mathematics and Statistics, University of Limerick, Ireland, and, Kevin Burke is a faculty member and Mark E. Burke an emeritus faculty member in the same department. B. Ross Barmish is a faculty member in ECE at Boston University and emeritus, also in ECE, at the University of Wisconsin. E-mail: Joseph.OBrien@ul.ie, Kevin.Burke@ul.ie, Mark.Burke@ul.ie, bob.barmish@gmail.com.} 
\thanks{The research of JDO'B was funded by Science Foundation Ireland Grant No. 16/IA/4470 and the initial work of BRB was funded by the MACSI Visitor Fund.}}

\begin{document}
	\maketitle
	\thispagestyle{empty}
	\pagestyle{empty}
	\vspace{-2in}
	\begin{abstract}
		For sequential betting games, Kelly's theory, aimed at maximization of the logarithmic growth of one's account value, involves optimization of the so-called betting fraction~$K$. In this Letter, we extend the classical formulation to allow for temporal correlation among bets. To demonstrate the potential of this new paradigm, for simplicity of exposition, we mainly address the case of a coin-flipping game with even-money payoff. To this end, we solve a problem with {\it memory depth}~$m$.  By this, we mean that the outcomes of coin flips are no longer assumed to be~i.i.d.~random variables. Instead, the probability of heads on flip~$k$ depends on previous flips~$k-1,k-2,...,k-m$. For the simplest case of~$n$ flips, with~$m = 1$, we obtain a closed form solution $K_n$ for the optimal betting fraction. This generalizes the classical result for the memoryless case. That is, instead of fraction~$K^* = 2p-1$ which pervades the literature for a coin with probability of heads~$p\geq 1/2$, our new fraction~$K_n$ depends on both~$n$ and the parameters associated with the temporal correlation. Generalizations of these results for~$m > 1$ and numerical simulations are also included. Finally, we indicate how the theory extends to time-varying feedback and alternative payoff distributions.
\end{abstract}

\section{Introduction\label{sec:intro}}

\IEEEPARstart{I}{n} Kelly's 1956 seminal paper~\cite{Kelly_1956}, the notion of {\it Expected Logarithmic Growth} (ELG) was introduced as the performance criterion for a memoryless repeated betting game. For a sequence of i.i.d.~gambles, for example a coin flip with the probability of heads being $p > 1/2$,  the theory leads to an {\it optimal betting fraction}~$K^*$, which, owing to its constant nature from bet to bet is viewed as a time-invariant feedback gain. That is, with~$V_k$ being the  account value after~$k$ plays, the optimal~$(k+1)$-th bet size is~$K^*V_k$, where, for classical coin flipping with even-money payoff, $K^* = 2p - 1$.

The ELG approach has resulted in a voluminous body of literature extending and applying the theory to other well-known gambling games such as Blackjack and sports betting considered in~\cite{Thorp_2006} and \cite{ haigh2000kelly},  asset management and stock trading as in  \cite{Cover_Thomas_2006, Luenberger_2013, rotando1992kelly,  Maslov98optimal, Hakansson_1971, lo2018growth, browne1996portfolio}; see also the extensive bibliography in~the~textbook~\cite{Maclean_Thorp_Ziemba_book}. Papers such as~\cite{Maclean_Thorp_Ziemba_2010, Maclean1992Growth, hsieh2016tooCons, RujeerapaiboonCEV18, Maclean2004security, busseti2016risk} have also covered related issues including asymptotics, problems related to aggressiveness of wagers and alternative risk metrics.

The main feature which differentiates this paper from existing work is our emphasis on the issue of {\it temporal correlation} among games. While it is standard to assume correlation among components of a multi-dimensional bet, for example, in modern portfolio theory~\cite{Luenberger_2013}, temporal correlation is an entirely different matter. Interestingly, although temporal effects are studied in the context of prediction for financial time-series, as in \cite{Campbell1993Correlation} and \cite{lewellen2002momentum}, this issue has received little attention in the context of bet sizing in Kelly's ELG framework; e.g., see \cite{Hirono2015temporal} where only one numerical example is considered.

Motivated by the fact that a bettor may gain an ``edge'' by taking advantage of temporal correlation, we generalize ELG theory to the case of a \textit{history-cognizant} coin where each outcome is no longer i.i.d.~but dependent on the previous~$m$~results. Our analysis of this new framework lays the groundwork for its use in financial applications with temporally correlated returns, in particular by relating the binary lattice model proposed in the sequel to stock price movements (up or down) over a sequence of time points. 
With this setting in mind, our primary analysis considers two-outcome, even-money random variables~$X_k \in \{-1,1\}$ with a time-invariant feedback gain governing the bet size which takes temporal autocorrelation into account. Although we also provide further extensions to accommodate $\ell$ possible outcomes given by~$X_k \in \{\mathrm{x}_1,\mathrm{x}_2,...,\mathrm{x}_{\ell}\}$ 
and time-varying feedback gains, our main focus is developing a new ELG theory in the presence of memory with arbitrary depth, $m \ge 1$.

The remainder of the paper is organized as follows: After formalizing the notion of {\it autocorrelated betting} in Section~\ref{Sec:autocorrelated}, we consider the {\it history-cognizant coin} in~Section~\ref{sec:histcoin}. Then, our main result and extensions are provided in~Section~\ref{Sec:mainresults}. This includes, for~$n$ bets and memory depth~$m = 1$, a closed-form solution for the optimal betting fraction~$K_n$, generalizing the classical~$K^* = 2p-1$ result, and steady state analysis; Section~\ref{Sec:proofs} is devoted to proof of this main theorem. Section~\ref{sec:deep} provides results for arbitrary memory depth ~$m \ge 1$ and model estimation, followed by numerical simulations and conclusions in the remaining two sections.

\section{Autocorrelated Kelly Betting}\label{Sec:autocorrelated}

\tikzstyle{block} = [draw, fill=orange!50, rectangle,
minimum height=3em, minimum width=5em]
\tikzstyle{sum} = [draw, fill=orange!50, circle, node distance=1cm]
\tikzstyle{input} = [coordinate]
\tikzstyle{output} = [coordinate]
\tikzstyle{pinstyle} = [pin distance=0.1cm, pin edge={to-,thin,white}]

\begin{figure}
	\centering	
	\begin{tikzpicture}[auto, node distance=2.6cm,>=latex',
	cross/.style = {circle, draw, fill=orange!50,
		append after command={\pgfextra{\let\LN\tikzlastnode}
			(\LN.north west) edge (\LN.south east)
			(\LN.south west) edge (\LN.north east)
	}}]
	\node [input, name=input] {};
	\node [cross, right of=input, node distance = 1.5cm] (multiplier) {};
	\node [sum, above of = multiplier, node distance = 1.3cm] (sum) {};
	\node [block, right of = multiplier, node distance = 2.5cm] (K) {K};
	\node [block, above of = K, pin = {[pinstyle]above:unit delay}] (delay) {$z^{-1}$};
	\node [output, right of=delay, label = right:{$V_k$}] (output) {};
	\draw [draw,->] (input) -- node[pos = 0.2] {$X_k$} (multiplier);
	\draw [draw, ->] (multiplier) -- node[pos = 0.85] {$+$}(sum);
	\draw [draw,->] (K) -- node[pos = 0.2] {} (multiplier);
	\draw [draw, ->] (sum) |- node[pos = 0.5] {$V_{k+1}$} (delay);
	\draw [->] (delay) -- node [name=y] {}(output);
	\draw [->] (y) |- (K);
	\draw [->] (y) |- node[above, pos = 0.97] {$+$} (sum);
	\end{tikzpicture}
	\caption{Feedback control configuration}
	\label{fig:feedback}
\end{figure}
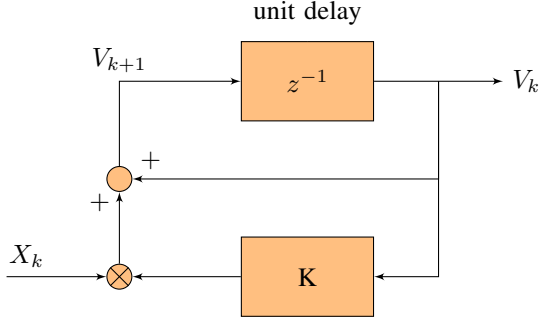

For the sake of self-containment, before introducing autocorrelation, we review Kelly's classical solution. Indeed, we start by considering a discrete-time even-money coin-flipping game with repeat i.i.d.~bets and initial account value~$V_0 > 0$. Letting $V_n$ denote the bettor's account value after $n$ plays, the classical \emph{Kelly} strategy is aimed at maximizing the expected value of the \emph{logarithm} of $V_n$ rather than simply its expected value. Letting $X_k \in \{-1,1\}$ be a random variable which represents the pay-off from the $k$th coin toss where, $X_k = 1$ corresponds to a head and $X_k = -1$ corresponds to a tail, the $(k+1)$-th bet is $KV_k$ with~$-1 \leq K \leq 1$. The quantity~$|K|$ is referred to as the betting fraction with $K<0$ representing a bet on tails rather than heads. Viewing~$V_k$ as a state, as noted for example in \cite{hsieh2016tooCons}, this defines a linear time-invariant feedback control~$u_k = KV_k$ for the nonlinear system as depicted in Figure~\ref{fig:feedback}, leading to the update
\begin{eqnarray*}
V_{k+1} & = & V_k + X_ku_k\\
& = & (1 + KX_k)V_k.
\end{eqnarray*}
We proceed to consider a game of $n$ bets with outcomes given by the sample path
$$
\mathbf{X} \doteq \left(X_0, X_1, X_2, \dots, X_{n-1}\right) \in \mathcal{X} \doteq \{-1,1\}^n.
$$
This being the case, the corresponding account value at terminal stage $n$, as a function of the pair~$(K,\mathbf{X})$, resulting from this sample path is given by
$$
V_n(K,\mathbf{X}) = V_0 \prod_{k=0}^{n-1}\left[1 + KX_k\right],
$$
and an optimal betting fraction is obtained by maximizing the Expected Logarithm Growth given by
$$
\text{ELG}(K) = \frac{1}{n} \mathbb{E}\left\{\log\left[\frac{V_n(K,\mathbf{X})}{V_0}\right]\right\}.
$$
Since~$\text{ELG}(K)$ above is independent of~$V_0$, in the sequel, without loss of generality, whenever convenient, it is assumed that~$V_0 = 1$. In the standard i.i.d.~setup, where $p$ is the probability of a head, the ELG is maximized at $K^* = 2p - 1$.

We now proceed to generalize the standard approach by assuming a probability distribution over~${\cal X}$ is available, and let $P_\mathbf{X}$ denote the probability of a sample path $\mathbf{X}$. This is a joint distribution over the components~$X_i$ of $\mathbf{X}$, and, at this high level of generality, this probability distribution is arbitrary. In the analysis to follow, we first provide a result for this general case which is abstract and then specialize to a scenario frequently encountered in practice. That is, we consider the case when the outcome of a given coin toss is correlated with the previous $m$ outcomes. We henceforth refer to~$m$ as the {\it memory depth} noting that a small value of $m$ means that the $k$th outcome is only related to the recent history. In this case, it is straightforward to see that the probability~$P_\mathbf{X}$ of a sequence $\mathbf{X}$ reduces to
\begin{align*}
P_\mathbf{X}
&= \prod_{k=0}^{n-1} \Pr(X_k | X_{k-1}, X_{k-2}, \ldots, X_{k-m})
\end{align*}
which is initialized by the $m$ events $X_{-m},X_{-m+1},\dots,X_{-1}$ prior the first outcome $X_0$ at stage $k = 0$.

\section{The History-Cognizant Coin\label{sec:histcoin}}
\begin{figure}
\centering
\includegraphics[width = 1.65in]{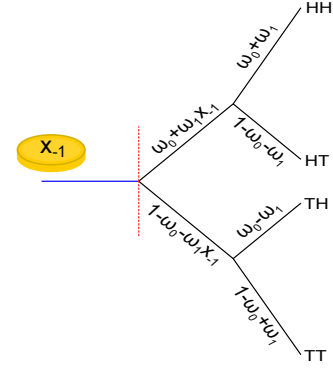}
\caption{Sample paths with $n = 2$ and $X_{-1} = x_{-1} \in \{-1,1\}$}
\label{fig:tree}
\end{figure}

Building on the above, we consider the simplest case of an autocorrelated bet: a coin whose current flip is affected by the previous flip. This is a coin with Markov memory, i.e., the probability of a head on the $k$th flip is
\begin{align*}
\Pr\left(X_k = 1 \, | \, X_{k-1}, X_{k-2},\dots\right) &= \Pr(X_k = 1  \, | \, X_{k-1}) \nonumber \\ &= \omega_0 + \omega_1 X_{k-1}
\end{align*}
with $\omega_0$ and $\omega_1$ assumed to be known; see Section \ref{sec:deep} for considerations of how these parameters may be estimated. For the parameterized linear function above, it is readily verified that the conditions~$~|\omega_1| < 0.5$, $|\omega_1| <~\omega_0 < 1 - |\omega_1|$ guarantee that $\Pr(X_k = 1  \, | \, X_{k-1}) \in (0,1)$. Now, via a straightforward calculation, these requirements reduce to
$$
\left|\omega_0 - \frac{1}{2}\right| + |\omega_1| < \frac{1}{2}
$$
which we recognize as describing the interior of  an~$\ell^1$ sphere, the so-called ``diamond'' with center~$(1/2,0)$ and radius~$(1/2)$. 
In this setting, we have memory depth~$m=1$, and we assume that we have observed one coin toss prior to betting, i.e., $X_{-1} = x_{-1} \in \{1,-1\}$. Figure~\ref{fig:tree} shows some illustrative sample paths, consistent with the formulation presented above. As mentioned in Section \ref{sec:intro}, this binary lattice can serve as a model for stock price movements over time categorized as ``up'' ($X_k = 1$) or ``down'' ($X_k = -1$).

\section{Main Result}\label{Sec:mainresults}
In this section, we provide our main result whose proof is relegated to the next section. The first part of the theorem below provides an abstract characterization of the optimal betting fraction $K_n$ in terms of the expected number of heads~$\mathbb{E}(H_n(\mathbf{X}))$ in the sample path $\mathbf{X}$ of length $n$, i.e., it holds for arbitrary sample path distributions $P_\mathbf{X}$ including and beyond those considered here. It also points the way to the second part of the theorem which addresses the case of a history-cognizant coin and makes use of notation
$$
p_0 \doteq \omega_0 + \omega_1 x_{-1}
$$
corresponding to the unconditional probability that~$X_0 = 1$, $$
{p}_\infty \doteq \frac{\omega_0 - \omega_1}{1 - 2\omega_1}
$$
which is later seen to be the steady state unconditional probability of heads, and
$$
\lambda_n \doteq \frac{1}{n}\left[\frac{1 - (2 \omega_1)^n}{1 - 2 \omega_1}\right]
$$
which satisfies the condition
$$
0 < \lambda_n < 1
$$
since ~$ |2\omega_1| < 1$ and tells us the relative weights of~$p_0$ and~$p_\infty$ in the optimal solution.

{\bf Theorem}: {\it For n flips of the history-cognizant coin with memory-depth $m = 1$ and conditional probability of heads given by $\Pr(X_k = 1  \, | \, X_{k-1}) \nonumber = \omega_0 + \omega_1 X_{k-1}$, the expected logarithmic growth~$\text{ELG}(K)$ is maximized by}
\begin{equation*}
K_n = 2\left\{\frac{\mathbb{E}(H_n(\mathbf{X}))}{n}\right\} - 1,
\label{eq:K_optim}
\end{equation*}
{\it where the expected value above is obtained as the convex combination}
$$
\frac{\mathbb{E}(H_n(\mathbf{X}))}{n} = \lambda_np_0
+ (1 - \lambda_n)p_\infty.
$$
\subsection{Special Cases, Generalizability and Remarks}

The remainder of this section focusses on finer details of our theory including its reduction to the classical Kelly formula for the memoryless case, results regarding the limiting values of the parameters used in the theorem, and generalizations of the theory beyond the simple case of time-invariant even-money two-outcome bets.

{\bf When Coin Flips are Independent}: For the special case with all payoffs~$X_k$ being i.i.d., we note that~$\mathbb{E}(H_n) = n p$ where $p = p_0 = p_\infty = \omega_0$ is the probability of a head. In this case,
$ K_n = K^* = 2p - 1,$
which is the classical result obtained in the absence of autocorrelation among bets as described in Section \ref{Sec:autocorrelated}.

{\bf Long Run Steady State Considerations}: Of general interest are the limiting values of the quantities $\mathbb{E}(H_n)$ and $K_n$ described in the theorem as~$n \rightarrow \infty$. The first point to note is that~$\lambda_n(\omega_1) \rightarrow 0$ which in turn implies that
$$
\lim_{n \rightarrow \infty}\frac{\mathbb{E}(H_n(\mathbf{X}))}{n} = p_\infty
$$
and immediately leads to optimal betting fraction
$$
K_\infty = 2p_\infty - 1.
$$
The interpretation of this limit is quite simple: If we are playing forever, the long-run probability of a head,~$p_\infty$, leading to~$K_\infty$ is the same betting fraction as that which one would obtain by ignoring correlation among the~$X_k$ and treating~$p_\infty$ as if it is the unconditional probability of a heads in the classical i.i.d.~case. On the other hand, if we are betting for a fixed time horizon~$n$, the difference between $\mathbb{E}(H_n) / n$ and $p_\infty$ is important. In particular, the optimal betting fraction~$K_n$ depends on~$n$ and the startup probability~$p_0 = \omega_0 + \omega_1 x_{-1}$, whereas $K_\infty$ does not. In practice, the importance of~$x_{-1}$ depends on the size of~$n$, and magnitude of autocorrelation coefficients $\omega_0$ and $\omega_1$.

{\bf Multiple Outcomes}: As indicated in Section \ref{sec:intro}, our theory may be modified to address more general scenarios. For example, consider the case where there are~$\ell$ possible outcomes~$\mathrm{x}_1,\mathrm{x}_2,...,\mathrm{x}_{\ell} \in (-1,\infty)$ for~$X_k$. Let~$P_{\mathbf{X}}$ be an arbitrary probability mass function over sample paths and let~$H_{n,i}(\mathbf{X})$ be a random variable denoting the number of times, in a path of length~$n$, that~$X_k = \mathrm{x}_i$ for~$ i = 1,2,...,\ell$. Then, 
using an argument quite similar to the one used in the proof of the theorem, we obtain
$$
\text{ELG}(K) = \sum_{i = 1}^{\ell} \frac{\mathbb{E}( H_{n,i}(\mathbf{X}))}{n}\log(1 + K \mathrm{x}_i),
$$
which is straightforward to maximize numerically since it is concave in $K$.

{\bf Time-Varying Feedback}: A second generalization which may also be considered, involves the use of time-varying feedback gains rather than the time-invariant $K$ synonymous with previous literature. The most straightforward extension in this direction is where, \textit{prior} to the start of the game, the bettor decides on a vector of betting fractions $$\mathbf{K} = \left[\tilde K_0, \tilde K_1, \ldots, \tilde K_{n-1}\right]^T.$$ Note that the tilde notation distinguishes these time-varying gains from the time-invariant $K_n$ of our main theorem. Thus, defining $p_k = \Pr(X_k = 1)$ to be the unconditional probability of a head on the $k$th coin toss,
\begin{align*}
\text{ELG}(\mathbf{K}) &= \frac{1}{n} \sum_{k = 0}^{n-1} \left\{p_k \log\left(1 + \tilde K_k\right) \right.\\[-0.2cm]
& \qquad\qquad\qquad + \left.(1-p_k) \log\left(1 - \tilde K_k\right) \right\},
\end{align*}
and this is maximized at $$\tilde K_{k} = 2p_k - 1.$$ Interestingly, it is straightforward to show that
$$K_n = \frac{1}{n}\sum_{k=0}^{n-1} \tilde K_{k} , $$
i.e., the time-invariant gain over $n$ bets is the average of the time-varying gains.

\section{Proof of the Theorem}\label{Sec:proofs}

This section can be skipped by those readers solely interested in the application aspects of this work. Indeed, to determine the optimal betting fraction, we maximize the Expected Logarithm Growth. For simplicity of notation, we suppress the dependence of~$H_n$ on the sample path~$\mathbf{X}$ and
calculate
\begin{eqnarray*}
\text{ELG}(K) &=& \frac{1}{n} \mathbb{E}\left\{\sum_{k = 0}^{n-1}\log(1 + K X_k)\right\} \\
&=& \frac{1}{n} \sum_{\mathbf{X} \in \mathcal{X}} P_{\mathbf{X}}\left\{\sum_{k = 0}^{n-1}\log(1 + K X_k)\right\} \\
&=&  \frac{1}{n} \sum_{\mathbf{X} \in \mathcal{X}} P_{\mathbf{X}} \left\{H_n \log(1 + K) \right. \\
& &\qquad\;\;\;\;\;\;\;\;\;\;\;\left. +(n - H_n) \log(1 - K)\right\} \\
&=& \frac{\mathbb{E}(H_n)}{n}\log(1+K)\\
& &\;\;\;\;\;\;\;\;\;\;\;+\left\{1 - \frac{\mathbb{E}(H_n)}{n}\right\}\log(1 - K).
\end{eqnarray*}
Now, noting that $K = K_n = 2\left\{\mathbb{E}(H_n) / n \right\} - 1$ is the unique point of zero derivative and that~$\text{ELG}(K)$ is a concave function, it follows that~$K_n$ is the unique maximizer.

It remains to derive an explicit formula for $\mathbb{E}(H_n)/n$ for the case of the history-cognizant coin. First, since the expected number of heads on the $k$th coin toss is $p_k = \Pr(X_k=1)$, the expected number of heads in $n$ coin tosses is
\begin{align*}
\mathbb{E}(H_n) = \sum_{k=0}^{n-1} p_k.
\end{align*}
Now, to obtain a formula for the sum above, beginning with conditional probability~$\Pr(X_k = 1  \, | \, X_{k-1}) = \omega_0 + \omega_1 X_{k-1}$, using the law of total expectation, we obtain a recursion
\begin{eqnarray*}
p_k 
& = & \mathbb{E}\left\{\Pr(X_k = 1|X_{k-1})\right\} \\
& = & \omega_0 + \omega_1 \mathbb{E}(X_{k-1})\\
& = &   2\omega_1 p_{k-1} + \omega_0 - \omega_1,
\end{eqnarray*}
where the last line follows since $\mathbb{E}(X_{k-1}) = 2 p_{k-1} - 1$. Initializing with~$p_0 = \omega_0 + \omega_1x_{-1}$,
the solution to this linear equation is, since~$|2\omega_1| < 1$,
\begin{align*}
p_k &= (2\omega_1)^k p_0 + \left[1-(2\omega_1)^k\right] p_\infty
\end{align*}
where $p_\infty = (\omega_0-\omega_1)/(1-2\omega_1)$. Thus,
\begin{align*}
\frac{\mathbb{E}(H_n)}{n} &= \frac{1}{n} \sum_{k=0}^{n-1} (2\omega_1)^k p_0 + \frac{1}{n} \sum_{k=0}^{n-1} \left[1-(2\omega_1)^k\right] p_{\infty}\\
&= \lambda_n p_0 + (1-\lambda_n) p_\infty. \;\;\square
\end{align*}

\section{Deeper Memory}\label{sec:deep}
We consider the general case of memory depth~$m > 1$ and show how analytic expressions for $\mathbb{E}(H_n)$ can be efficiently obtained. Indeed, beginning with parameterization of the conditional probability of heads
$$
\Pr(X_k = 1|X_{k-1},X_{k-2},...,X_{k-m}) = \omega_0 + \sum_{i = 1}^m\omega_i X_{k-i}
$$
with assumed initial conditions
$$
X_{-i} = x_{-i} \;\mbox{for}\; i = 1,2,...m.
$$
To ensure that~$\Pr(X_k = 1|X_{k-1},X_{k-2},...,X_{k-m}) \in (0,1)$ we assume the~$\omega$ parameters to lie in the ``hyperdiamond''
$$
\left|\omega_0 - \frac{1}{2}\right| + \sum_{i=1}^m |\omega_i| < \frac{1}{2}.
$$
Then, the unconditional probability~$p_k$ that~$X_k = 1$ is
\begin{align*}
p_k &= \mathbb{E}\{\Pr(X_k = 1|X_{k-1},X_{k-2},\dots,X_{k-m})\}, \\
&= \omega_0 + \sum_{i = 1}^m\omega_i \mathbb{E}(X_{k-i})
\end{align*}
for~$k = 0,1,\dots,n-1$. Substituting $\mathbb{E}(X_{k-i}) = 2p_{k-i}-1,$
above and taking note of the ``induced'' initial conditions

$$
p_{-i} = \frac{(x_{-i}+1)}{2}\; \mbox{for}\; i = 1,2,\dots,m,
$$
we arrive at the recursion
\begin{align*}
p_k &= \omega_0 - \sum_{i=1}^m \omega_i + 2 \sum_{i=1}^m \omega_i p_{k-i}
\end{align*}
which holds for~$k = 0,1,...,n-1$, and from which $E(H_n)$ and hence $K_n$ may be calculated. To illustrate a specific case, for memory depth~$m= 3$ and~$n = 2$ flips, we find that
\begin{align*}
\mathbb{E}(H_2) &= x_{-1}(2\omega_1^2 + \omega_1 + \omega_2) + x_{-2}(2\omega_1\omega_2 +\omega_2 +\omega_3)\\
&\quad + x_{-3}( 2\omega_1\omega_3 + \omega_3) + 2\omega_0\omega_1 + 2\omega_0 - \omega_1.
\end{align*}

{\bf State-Space Formulation}: As an alternative to the above, which may perhaps prove useful in future research, we consider a standard state-space realization of the ``delay system'' to represent the scalar recursion for $p_k$.  That is, by introducing the~$m$-dimensional state vector which is given by~$v_k = [p_{k-m+1}, p_{k-m + 2},\ldots,p_{k}]^T$, we readily obtain a classical companion form realization ~$v_{k+1} = Av_k + b u_k$ with triple~$(A,b,c)$,  input~$u(k) \equiv 1$ and output being~$p_k$. To illustrate, for memory depth~$m = 3$, we obtain
$$
A = \left[
\begin{array}{ccc}
0 & 1 & 0 \\
0 & 0 & 1 \\
2\omega_3 & 2\omega_2 & 2\omega_1 \\
\end{array}
\right];\;\;
b = \left[
\begin{array}{c}
0 \\
0 \\
\omega_0 - \omega_1 - \omega_2 - \omega_3\\
\end{array}
\right];
$$
and~$c = [0\;0\;1]$ with solution of the recursion given by
\begin{align*}
p_k & =  c\left(A^k v_0 + \sum_{i=0}^{k-1}A^{k-1-i}b \right) \\
& =  c\left(A^k v_0 + (I-A)^{-1}(I-A^k)b\right)
\end{align*}
and the matrix~$I-A$ is guaranteed to be invertible since $\text{det}(I- A) = 1 - 2\sum_{i=1}^m \omega_i$
must be non-zero due to the hyperdiamond constraint
on the $\omega_i$. In addition, from Gerschgorin's circle theorem~\cite{golub2012matrix} and the hyperdiamond constraint, each eigenvalue of~$A$ has magnitude less than 1, and so~$A^k \rightarrow 0$ as~$k \rightarrow \infty$. Using this fact, beginning with $p_k$ above, this leads to the generalization of our steady-state unconditional probability
$$ p_{\infty} = c(I - A)^{-1}b. $$
Now recognizing that this corresponds to the transfer function $H(z)$ for the triple $(A,b,c)$ evaluated at $z = 1$, we immediately arrive at
$$p_\infty = \frac{\omega_0 - (\omega_1 + \dots + \omega_m)}{1 - 2(\omega_1 + \dots + \omega_m)}.$$


{\bf Estimation}: In practice, prior to betting, it is necessary to obtain values for the $\omega_i$ parameters. First, define the ``response'' variable $Y_k = (X_k+1)/2$
such that
$$
\mathbb{E}(Y_k  \, | \, X_{k-1}, X_{k-2},\dots,X_{k-m}) = \omega_0 + \sum_{i = 1}^m\omega_1 X_{k-i}.
$$
Then, having observed data $x_{-\ell},\dots,x_{-1}$, we compute the $(\ell-m) \times 1$ response vector
$$
y \doteq \frac{1}{2}\left\{ [x_{-\ell+m},\, x_{-\ell+m+1},\, \ldots,\, x_{-1}]^T + 1\right\},
$$
and minimize the residual sum of squares
$$
\text{RSS}(\omega) = \sum_{k=-\ell+m}^{-1} \left(y_k - \omega_0  - \sum_{i = 1}^{m}\omega_i x_{k-i}\right)^2,
$$
with respect to $\omega = [\omega_0, \omega_1, \ldots, \omega_m]^T.$ Whereas classical estimation theory leads to the \emph{least squares} solution
$$
\hat \omega = \underset{\omega}{\operatorname{argmin}}~\text{RSS}(\omega) = (X^TX)^{-1}X^Ty
$$
with $X$ being the $(\ell-m) \times (m+1)$ matrix whose $i$th row is given by~$
[1,\, x_{i-\ell+m-2},\, x_{i-\ell+m-3},\, \ldots,\,  x_{i-\ell-1}]$, enforcement of the hyperdiamond constraints leads to a positive-definite convex program to be solved.

\begin{figure*}[htb]
\centering
\includegraphics[width = \textwidth, trim = {0 0 0 0}]{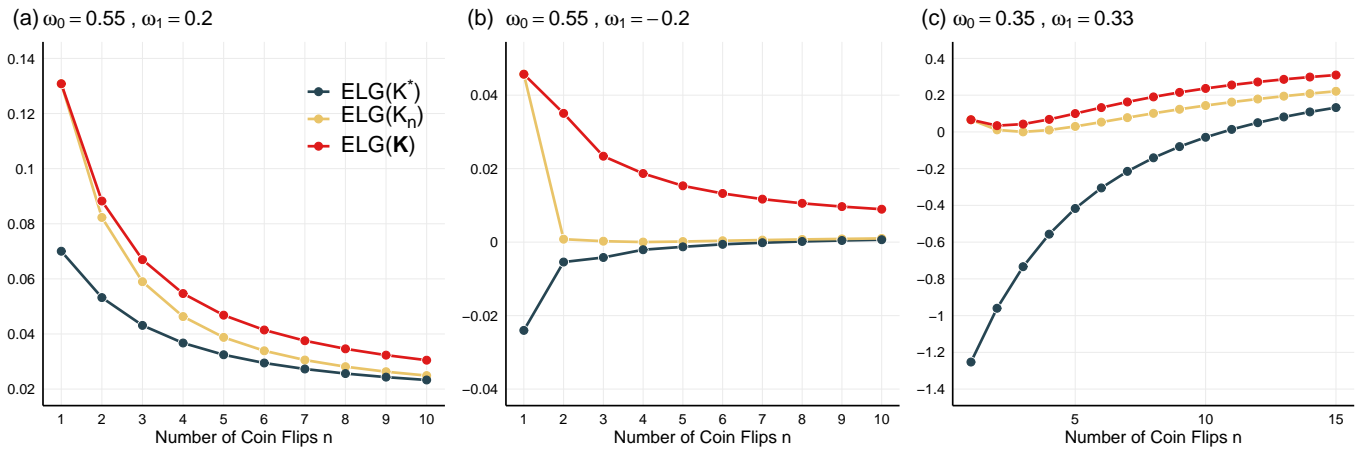}
\caption{Simulation demonstrating the ELG of account value for the three bettors with $X_{-1} = 1$, and different combinations of $\omega_0$ and $\omega_1$.}
\label{fig:simulation}
\end{figure*}

\section{Numerical Simulation\label{Sec:numerical}}

We present results for simulations with returns driven by a process with $\Pr(X_k = 1 \, | \, X_{k-1}) = \omega_0 + \omega_1X_{k-1}$, comparing ELG performance for the classical Kelly $K^*$-bettor, the $K_n$-bettor, and the $\mathbf{K}$-bettor. Accordingly, from Section \ref{Sec:mainresults}, the classical $K^*$-bettor, disregarding temporal correlation, works with the probability of heads being $p_\infty$ and uses the time-invariant betting fraction~$K^* = 2p_\infty -1$.
On the other hand, assuming $\omega_0$ and~$\omega_1$ are perfectly estimated, the $K_n$-bettor exploits temporal correlation and uses fraction $K_n = 2\{\lambda_n p_0 + (1-\lambda_n)p_\infty\} - 1$, while the $\mathbf{K}$-bettor, who is also aware of the correlation, uses separate $K$ values for each stage via $\tilde K_k = 2 p_k - 1$.

To provide a flavor of our findings, we first consider the following scenario: We initialize the game by supposing that the prior event was~$X_{-1} = 1$, and take~$\omega_0 = 0.55$, $\omega_1 = 0.20$. In this specific situation, a straightforward calculation using our theory directly leads to values of~$p_0 = 0.75, p_\infty = 0.583$, and $\lambda_n \approx 0.556(1 - 0.4^n)/n$; hence, $K^* \approx 0.167$ and $K_n \approx 0.556 (1 - 0.4^n)/n + 0.167$
from which it is immediately clear that $K_n > K^*$ (but tends to $K^*$ in the long run per Section \ref{Sec:mainresults}). In this setting, due to the positive correlation ($\omega_1 > 0$), and fact that $X_{-1} = 1$, earlier coin tosses are more likely to be heads than later ones. Therefore, when betting for a fixed time horizon $n$, the $K_n$-bettor takes advantage of the temporal correlation by placing larger bets than those suggested by the correlation-ignoring $K^*$ value. However, although accounting for autocorrelation, $K_n$ is time-invariant. For $n=2$, $K_n = 0.4$, whereas $\mathbf{K} = (0.5, 0.3)$ (and recall from Section \ref{Sec:mainresults} that $K_n$ is the average of the elements of $\mathbf{K}$). Thus, the $\mathbf{K}$-bettor, recognizing that $X_0$ is most likely to be a head, bets more heavily on the first bet than on the second. Of course, both the $K_n$- and  $\mathbf{K}$-bettors bet more heavily than the $K^*$-bettor, and, indeed, for the $n = 2$ case, we find that $\text{ELG}(K^*) \approx 0.053$, $\text{ELG}(K_n) \approx 0.082$, and $\text{ELG}(\mathbf{K}) \approx  0.088$.

Beyond $n=2$, Figure~\ref{fig:simulation}(a) shows ELG values for a range of $n$ over which the $\mathbf{K}$-bettor outperforms the $K_n$-bettor who in turn outperforms the $K^*$-bettor. This scenario is analogous to one which arises for a financial asset on an upward  trend (since $p_k \ge p_\infty \approx  0.5833$). In such a setting, the majority of strategies will do well, e.g., all three here have positive ELG, but, importantly, incorporating temporal correlation boosts performance. Figure~\ref{fig:simulation}(b) displays the results for a similar simulation but with $\omega_1 = -0.2$. As with the first scenario, this represents long-run upward trend since $p_\infty \approx  0.5357$, but the negative autocorrelation means that the process fluctuates more; in particular, $X_0$ is most likely to be a tail since $X_{-1}$ was a head. In this setting, neither the $K^*$- nor the $K_n$-bettors do very well, albeit the latter at least has non-negative ELG, whereas the $\mathbf{K}$-bettor has significantly improved performance. To see why this is, consider the $n=2$ case where $K^* \approx  0.071$, $K_n = -0.04$, and $\mathbf{K} = (-0.3, 0.22)$. Thus, the $\mathbf{K}$-bettor makes use of the fluctuation by betting on tails first and on heads second but bets less heavily in the second due to the increased uncertainty; the $K_n$-bettor averages over these fluctuations, slightly favouring tails but ultimately betting very little, whereas the $K^*$-bettor wrongly favours heads.

Lastly, in Figure \ref{fig:simulation}(c) we briefly consider another interesting scenario shown which corresponds to $X_{-1} = 1$, $\omega _0 = 0.35$ and~$\omega_1 = 0.33$. The important feature in this case is that it corresponds to a long-run downward trend (with $p_\infty = 0.058$) but where the positive temporal correlation and the fact that $X_{-1} = 1$ mean that $X_0 = 1$ is most likely. In this case, the $K^*$-bettor suboptimally takes a heavy short position (i.e., bets on tails) with $K^* \approx -0.883$ leading to large negative ELG value early on (only becoming positive for $n > 10$). In contrast, both $K_n$- and $\mathbf{K}$-bettors achieve positive growth over all $n$. Note that, for all three scenarios, and over all $n$ values, $\text{ELG}(K^*) \le \text{ELG}(K_n) \le \text{ELG}(\mathbf{K})$ which is consistent with our exploitation of temporal correlation to obtain improved Kelly-type betting strategies which have not been considered in the existing literature.

\section{Conclusion}
In this paper, we formulated a class of Kelly optimal ELG problems which account for temporal correlation over the sequence of gambles. In the main theorem, for memory depth~$m =1$ and~$n$ flips, a closed form solution for the optimal betting fraction~$K_n$ was obtained. The paper also includes analysis for the case when~$n \rightarrow \infty$ and solutions for deeper memory~$m > 1$ which can be obtained by either propagation of the recursive formula for $p_k$ or use of the state space realization for the associated delay system. While our primary focus has been the development of ELG theory in the presence of autocorrelation, we have also provided extensions to multiple-payoffs and time-varying feedback gains. Numerical simulations which included comparison with classical Kelly betting results on games with temporal correlation were also shown in order to demonstrate the potential advantages offered by our framework.

In future work, we envision our theory to be especially applicable to scenarios in which an investor wishes to incorporate temporal correlations into algorithmic trading strategies over the course of time. In particular, the conceptual framework introduced within this Letter has potential to provide the base upon which multiple extensions beyond those proposed above can be built; one such avenue is the so-called {\it portfolio case} with correlation both temporally and across components. In this case it is felt that concave programming will play an important role in computation; e.g., see \cite{boyd2004convex}. Finally, a further direction of research involves a study of ELG performance as a function of betting frequency in the context of the temporal correlation framework introduced here; e.g., see\cite{Wu2015adaptive} for initial work along these lines and \cite{hsieh2018frequency} for analysis of the memoryless case.

\bibliographystyle{unsrt}

\end{document}